# A discussion on numerical shock stability of unstructured finite volume method: Riemann solvers and limiters

Fan Zhang[1, a *], Zhichao Yuan[2,b], and Jun Liu[2,c]

[1]Department of Engineering Mechanics, Dalian University of Technology, Dalian, China

[2]School of Aeronautics and Astronautics, Dalian University of Technology, Dalian, China

[a]a04051127@mail.dlut.edu.cn, [b]zhichaoy@126.com, [c]liujun65@dlut.edu.cn

**Abstract:** Numerical shock instability is a complexity which may occur in supersonic simulations. Riemann solver is usually the crucial factor that affects both the computation accuracy and numerical shock stability. In this paper, several classical Riemann solvers are discussed, and the intrinsic mechanism of shock instability is especially concerned. It can be found that the momentum perturbation traversing shock wave is a major reason that invokes instability. Furthermore, slope limiters used to depress oscillation across shock wave is also a key factor for computation stability. Several slope limiters can cause significant numerical errors near shock waves, and make the computation fail to converge. Extra dissipation of Riemann solvers and slope limiters can be helpful to eliminate instability, but reduces the computation accuracy. Therefore, to properly introduce numerical dissipation is critical for numerical computations. Here, pressure based shock indicator is used to show the position of shock wave and tunes the numerical dissipation. Overall, the presented methods are showing satisfactory results in both the accuracy and stability.

## Introduction[1]

Numerical simulations on strong shock waves have proved to be very challenging because of the appearance of shock anomalies, such as carbuncle phenomenon [1], odd-even decoupling [2], and so on. Many researches have contributed to discussions on the mechanism of shock anomalies and/or improvements on numerical flux schemes. However, debates and investigations are still carrying on these issues. For example, it has been conjectured that the carbuncle phenomenon is a physical instability. Elling [3] had triggered carbuncle-like numerical results by resembling a vortex sheet in the upstream flow of shock wave. Dumbser, *et al.* [4] suggested that the numerical carbuncle phenomenon is due to an unconditional instability of the underlying mean flow, and infinitely small errors will grow exponentially in time if they belong to the unstable modes.

Some of the aforementioned anomalies, e.g. carbuncle phenomenon, usually investigated by using first-order spatial discretization. Slope limiter is a crucial ingredient to keep numerical stability of second-order schemes. However, slope limiter does not guarantee stability due to the complexity of physical phenomena and the discretization.

In this paper, numerical behavior of several typical upwind flux schemes and slope limiters are investigated. The results are discussed and the improvements are compared with the traditional methods.

## Discretization of Hypersonic Flow

Hypersonic flow around a cylinder is a typical situation that causes shock instability, i.e. carbuncle phenomenon. This phenomenon has a possible physical background which was introduced by Kalkhoran, *et al.* [5]. In their experiment, a weak perturbation of vortex sheet travelling downstream to a normal shock wave, and a leading oblique shock is produced.

Here, three discretization, e.g. quadrilateral, regular triangle, and irregular triangle, are shown in Fig.1. The sparse grids in the figures are schematic, and the grids applied in the simulated are refined.

---

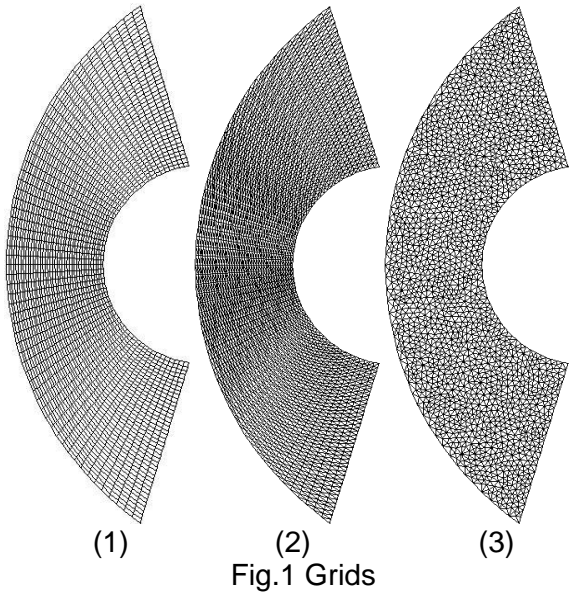
(1)      (2)      (3)
Fig.1 Grids

**Primary Numerical Results and discussion**

Mach number of this case is set as $Ma=8$, and inviscid Euler equations are used as the governing equations. Cell-centered Finite Volume (FV) method in [6] is applied for spatial solution.

The first-order results of Roe scheme [7] and van Leer scheme [8] on quadrilateral grid are shown in Fig. 2, where a typical carbuncle is produced by Roe scheme, yet the van Leer scheme shows a clear and stable bow shock.

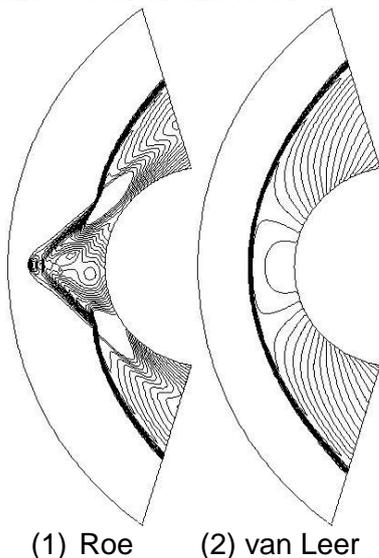
(1) Roe      (2) van Leer
Fig.2 Numerical results on quadrilateral gird

Many upwind schemes, e.g. van Leer, AUSM$^+$ [9], are usually showing stable results on quadrilateral grid. However, it is difficult to build high quality symmetrical grids in real-life applications. A typical situation is the triangular grids, e.g. grid 2 or 3 in Fig 1.

On these types of grids, Roe scheme is not the only one which produces anomalies. The AUSM$^+$ scheme, and SLAU scheme [10] which also belongs to the AUSM family, had shown less satisfactory results, which are shown in Fig.3.

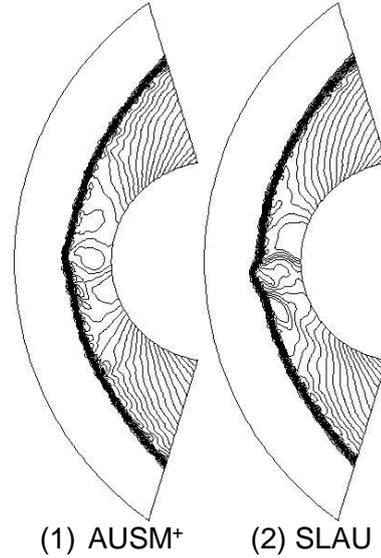
(1) AUSM$^+$      (2) SLAU
Fig.3 Numerical results on irregular triangular gird

The instability occurred on certain discretization is probably due to unsymmetrical flux (and perturbation) caused by the grid. The problem can be observed specifically by using the grid 2 in Fig.1, which is produced by diagonally cutting the quadrilateral cell of grid 1. The result of SLAU scheme on this regular triangular grid is shown in Fig. 4.

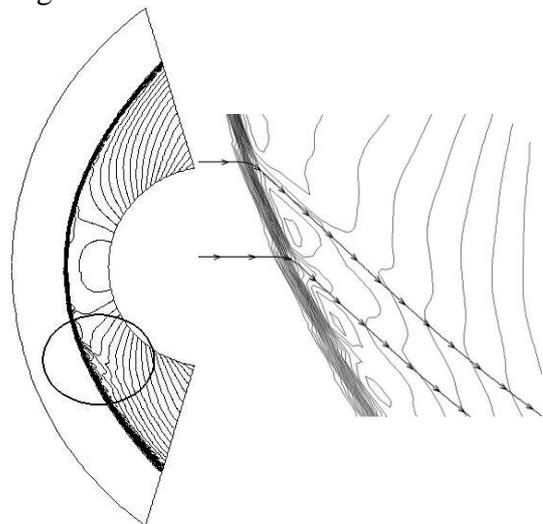
Fig.4 Numerical results of SLAU scheme on regular triangular gird

In this case, SLAU scheme shows anomaly behind the bow shock. Detailed discussion was introduced in [9]. In general, although the anomaly is not as significant as the carbuncle phenomenon, it is produced while the grid line and stream line are coincided and obliquely cutting the shock wave, and thus the left and right variables of the grid line are significantly different, which travels downstream if the upwind scheme keeps shear wave and contact discontinuity perfectly. The explanation is somehow fits the result of Elling [3] and Kalkhoran, *et al*. [5], which induce shock instability by resembling a vortex in upstream.

**Numerical Results of new schemes or different discretization**

In [11], a pressure-based function was introduced to suppress the instability, which is

$$\begin{cases} \omega = \min_{l \in k(i) \cup k(j)} (f_l) \\ f_l = \min\left(\dfrac{p^+}{p^-}, \dfrac{p^-}{p^+}\right)_l^3 \end{cases} \quad (1)$$

where the parameter is shown in Fig. 5.

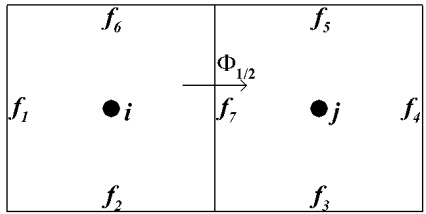

Fig.5 Illustration of the definition of the weights in hybrid scheme

This function has also been applied to improve the stability of TV scheme [12,13]. Numerical results of improved SLAU scheme and TV scheme are shown in Fig. 6.

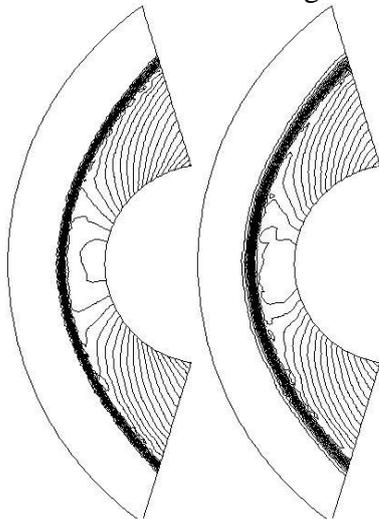

(1) SLAU-Hybrid (2) TV-hybrid
Fig.6 Numerical results of hybrid schemes

The essential idea here is to increase dissipation while pressure discontinuity exists. However, this effective *ad hoc* strategy is not further discussed. It should be mentioned that the momentum perturbation in numerical fluxes has been deemed to be the main source of shock instability, and thus the hybrid scheme only needs to increase the dissipation for momentum flux.

Furthermore, more numerical results are introduced to prove the travelling perturbation plays a key role in the shock stability. In Fig. 7, vertex-centered scheme is used to discretize the flow field, with using the same grid in Fig. 1.

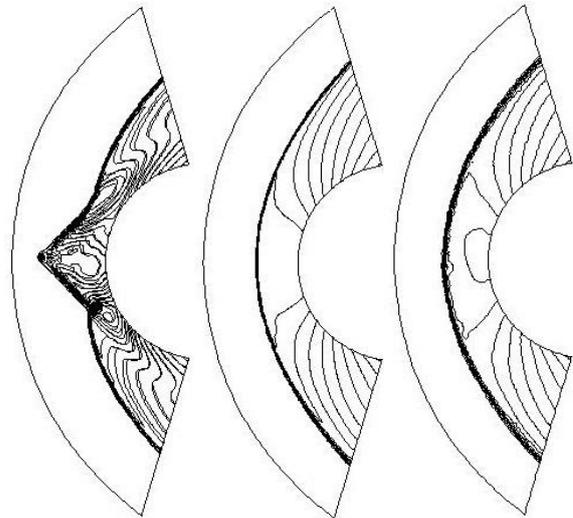

(1) grid1    (2) grid 2    (3) grid 3
Fig.7 Numerical results of vertex-centered FVM (Roe scheme)

Still, Roe scheme shows carbuncle on quadrilateral grids, since the vertex-centered scheme also have quadrilateral control volume on this grid. Whereas, on triangular grids, control volume is hexagon for the vertex-centered FVM. Correspondingly, carbuncle is eliminated on triangular grids. As a matter of fact, these results are reasonable because the perturbation travelling path is restricted on hexagon discretization. Therefore, although the perturbation will be produced and the upwind scheme (Roe) does not suppress the perturbation, the numerical results are still stable on triangular grids.

**Numerical Results of Second-Order Schemes**

Aforementioned first-order results give various conditions that produce carbuncle. Usually, for second-order schemes carbuncle phenomenon is alleviated, although not eliminated. Since the second-order scheme reduces discontinuities between cells, the corresponding perturbation will be less significant. However, anomalies still exist.

Fig. 8 present three results of typical unstructured slope limiters, including Venkatakrishnan limiter [14], Barth limiter [15] and MLP limiter [16].

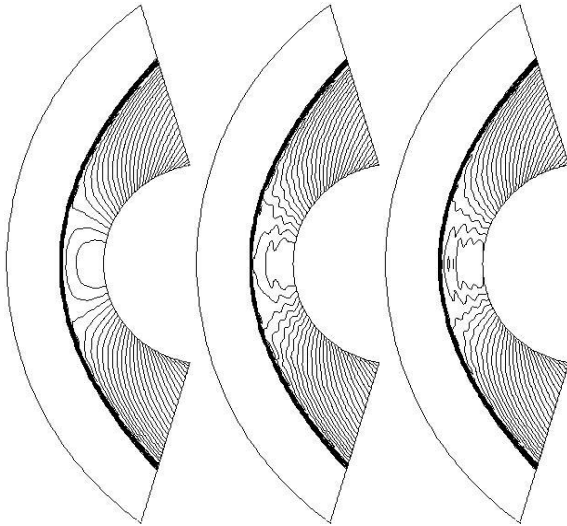

(1) Venkatakrishnan (2) Barth (3) MLP
Fig.8 Second-order schemes on quadrilateral grid

Here, Venkatakrishnan limiter and MLP limiter have been given the parameter $K=1$, and convective flux is calculated by van Leer scheme, which is deemed to be stable for simulating shock wave. Therefore, the oscillations here should be induced by slope limier.

The improved version of MLP had been introduced in [17]. A strict-MLP condition and a weak-MLP condition are proposed in this paper, and these two conditions are combined by a pressure-based function similar with Eq.1.

The idea of the strict-MLP condition is to restrict the reconstructed interface variables to a medium value that is between the left and right cell-centre values, and thus the reconstruction will not cause new extrema. The condition is schematically shown in Fig. 9.

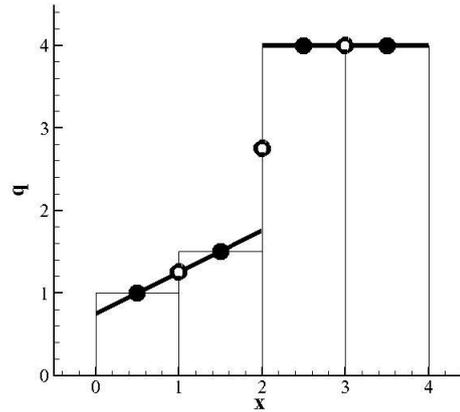

Fig.9 Strict-MLP condition in 1D

The numerical results of MLP-pw limiter are shown in Fig.10. With setting the parameter $K=1$, or $K=10$, the results are both stable.

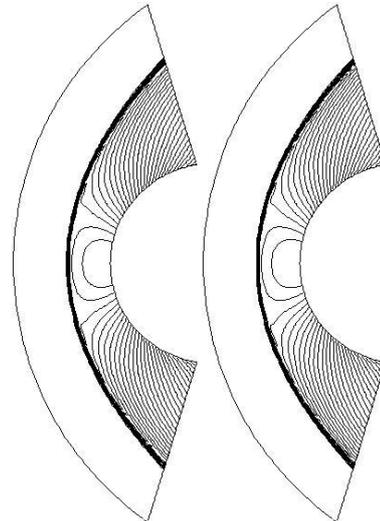

(1) $K=1$      (2) $K=10$
Fig.10 Numerical results of MLP-pw limiter on quadrilateral grid

The convergence of several traditional limiter and the improved limiter are shown in Fig.11.

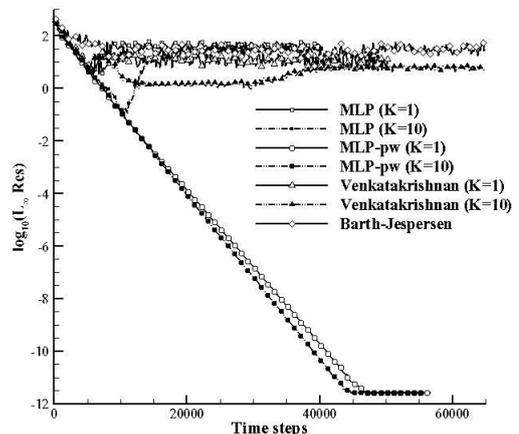

Fig.9 Computational residual on quadrilateral grid

It can be found that the improved limiter removes the oscillations in the flow field, and thus the convergence is significantly improved.

**Conclusions**

Numerical stability are the results of upwind flux function and slope limiter. For first-order schemes, slope limiter is not applied and thus the effects of upwind flux scheme can be specifically observed. The instability is caused by the co-effect of grid discretization, flow condition and dissipation of the numerical scheme. Increase the dissipation is an effective way, and different discretization which reduces the travelling path of perturbation can also eliminate the instability.

For second-order scheme which needs slope limiter to suppress Gibbs phenomenon, upwind flux scheme is not the only reason causes instability. Especially, instability exists even a stable upwind scheme has been used. In this paper, the results of an improved limiter are shown, in which the instability is eliminated by properly introducing dissipation near shock wave.